\documentclass[smallextended,envcountsect]{svjour3}
\smartqed
\usepackage{graphicx}
\journalname{JOTA}

\usepackage{amsfonts}
\usepackage{amssymb}
\usepackage{graphicx}
\usepackage{cancel}
\usepackage{tikz}

\def\emp{\emptyset}
\def\conv{{\rm conv}\,}

\def\dom{{\rm dom}\,}
\def\span{{\rm span}\,}
\def\epi{{\rm epi\,}}

\def\N{{\cal N}}

\def\O{{\cal O}}

\def\sub{\partial}
\def\B{\mathbb B}

\def\ox{\overline{x}}
\def\oy{\overline{y}}
\def\oz{\overline{z}}

\def\disp{\displaystyle}

\def\tto{\;{\lower 1pt\hbox{$\rightarrow$}}\kern-10pt
\hbox{\raise 2pt\hbox{$\rightarrow$}}\;}

\def\Bar{\overline}
\def\ra{\rangle}
\def\la{\langle}

\def\B{I\!\!B}

\def\R{\mathbb{R}}
\def\N{I\!\!N}
\def\ox{\bar{x}}
\def\oy{\bar{y}}
\def\oz{\bar{z}}
\def\ov{\bar{v}}
\def\ow{\bar{w}}

\def\oq{\bar{q}}
\def\co{{\rm co}}
\def\cone{{\rm cone}}

\def\int{\mbox{\rm int}\,}
\def\gph{\mbox{\rm gph}\,}
\def\epi{\mbox{\rm epi}\,}
\def\dim{\mbox{\rm dim}\,}
\def\dom{\mbox{\rm dom}\,}

\def\ker{\mbox{\rm ker}\,}

\def\aff{\mbox{\rm aff}\,}

\def\conv{\mbox{\rm conv}}

\def\b{\hfill\Box}

\def\O{\Omega}
\def\ph{\varphi}
\def\emp{\emptyset}
\def\st{\stackrel}
\def\oR{\Bar{\R}}
\def\lm{\lambda}
\def\gg{\gamma}
\def\dd{\delta}
\def\al{\alpha}
\def\Th{\Theta}
\def \N{I\!\!N}
\def\th{\theta}
\def\vt{\vartheta}
\def\beq{\begin{equation}}
\def\eeq{\end{equation}}

\begin{document}
\title{Second-Order Analysis of Piecewise Linear Functions with Applications to Optimization and Stability}
\author{B. S. Mordukhovich \and M. E. Sarabi}

\institute{Boris S. Mordukhovich  \at
Wayne State University \\
Detroit, Michigan\\
\email{boris@math.wayne.edu}
\and
M. Ebrahim Sarabi,  Corresponding author  \at
Wayne State University \\
Detroit, Michigan\\
\email{ebrahim.sarabi@wayne.edu}}
\dedication{\rm\small{Communicated by Aram Arutyunov }}
\date{Received: date / Accepted: date}

\maketitle

\begin{abstract}
This paper is devoted to second-order variational analysis of a rather broad class of extended-real-valued piecewise liner functions and their applications to various issues of optimization and stability. Based on our recent explicit calculations of the second-order subdifferential for such functions, we establish relationships between nondegeneracy and second-order qualification for fully amenable compositions involving piecewise linear functions.
We then provide a second-order characterization of full stable local minimizers  in composite optimization and constrained minimax problems.
\end{abstract}
\keywords{variational analysis and optimization \and piecewise linear functions \and second-order subdifferentials \and nondegeneracy \and full stability of local minimizers}
\subclass{49J52 \and 90C30 \and 90C31}

\section{Introduction}

This paper concerns developing second-order generalized differential theory of variational analysis and its applications to problems of nondifferentiable optimization as well as to several notions of stability in parametric optimization and associated variational systems. Professor Vladimir Demyanov made very significant contributions to these areas (see, e.g., \cite{dem} and the references therein) that have been well recognized by the optimization community.

In this paper we mainly address variational theory and applications of the class of {\em convex piecewise linear} (CPWL) extended-real-valued functions \cite{rw} playing an important role in many aspects of variational analysis and optimization. Having in hands recently obtained \cite{ms15} explicit calculations of the {\em second-order subdifferentials} (or generalized Hessians) of such function in the sense of \cite{m92}, we present here some of their applications to second-order variational analysis and parametric optimization. Proceeding in this direction requires us to deal not only with CPWL functions per se but mainly with {\em fully amenable compositions} involving such functions, which play an underlying role in many aspects of variational analysis, optimization, and stability.

The first issue studied in this paper is to clarify relationships between the {\em exact} (equality-type) {\em second-order chain rules} derived recently in \cite{mr} and \cite{mnn} for fully amenable compositions under different qualification conditions. Using the second-order calculations from \cite{ms15} allows us to show that these two conditions are equivalent in a certain precise sense in the case of convex piecewise linear outer functions. This leads us to a deeper understanding of second-order variational calculus and its applications.

The next topic we address in this paper is {\em full stability} (in the sense introduced by Levy, Poliquin and Rockafellar \cite{lpr} in the general extended-real-valued framework of unconstrained optimization) of local minimizers for a rather broad class of {\em composite optimization} problems governed by fully amenable compositions with CPWL functions. Employing again the second-order calculations from \cite{ms15} and the second-order subdifferential sum and chain rules leads us to deriving complete characterizations of fully stable local minimizers for such composite problems expressed entirely in terms of their initial data via the appropriate {\em composite SSOSC} (strong second-order sufficient condition) under a certain {\em partial nondegeneracy}. An effective implementation of this result is given for the case of {\em constrained minimax problems}.

The rest of the paper is organized as follows. In Section~2, we present some preliminaries and recall, for the reader's convenience, basic definitions and results from \cite{ms15} required in the sequel. It makes this paper to be  fully {\em self-contained}, and the reader needs to consult the the related paper \cite{ms15} only for the proofs of the underlying results used here for applications.\vspace*{-0.02in}

The major result of Section~3 establishes a certain equivalence between qualification conditions used in \cite{mr,mnn} for deriving by different approaches the aforementioned exact second-order chain rule for fully amenable compositions involving CPWL functions. Being important for its own sake, the key ingredient of this result (together with the explicit calculation of the second-order subdifferential of CPWL functions) is the proof of the so-called {\em ${\cal C}^\infty$-reducibility} of CPWL functions via {\em linear} transformations that are used then in the formulation of partial nondegeneracy. As a by-product of the obtained equivalence, we completely clarify the essence of the powerful second-order chain rule that is largely employed in the subsequent material.\vspace*{-0.02in}

Section~4 presents the explicit composite SSOSC characterization of {\em fully stable} local minimizers in the partially nondegenerate {\em composite framework} of optimization involving CPWL functions. In Section~5, we effectively apply the general composite result to characterizing full stability of local solutions {\em minimax} problems with {\em polyhedral constraints}.\vspace*{-0.02in}

Throughout the paper we use the standard notation of variational analysis; see \cite{rw,m06}. For brevity, a number of acronyms is applied in the next. Besides those mention above, they include: SOCQ (second-order qualification condition), LICQ (linear independence constraint qualification),  AICQ (affine independence constraint qualification), MPPCs (mathematical programs with polyhedral constraints), SOCPs (second-order cone programs), NLPs (nonlinear programs), and ENLPs (extended nonlinear programs).\vspace*{-0.2in}

\section{Basic Definitions and Preliminaries}\label{prel}

Recalling first the constructions of generalized differentiation used below, we begin with the basic definition of generalized normals to arbitrary sets, the only one normal cone construction employed in the paper. Given $\O\subset\R^n$ with $\ox\in\O$, the {\em normal cone} to $\O$ at $\ox$ (known also as the limiting, basic or Mordukhovich normal cone) is defined by
\begin{equation}\label{2.2}
N(\ox;\O):=\Big\{v\in\R^n\colon\exists x_k\to\ox,v_k\to x\;\mbox{s.t.}\;x_k\in\O,\disp\limsup_{x\st{\O}\to x_k}\frac{\la v_k,x-x_k\ra}{\|x-x_k\|}\le 0\Big\},
\end{equation}
where $k\in\N:=\{1,2,\ldots\}$. It is well known that, despite the intrinsic nonconvexity of (\ref{2.2}) for nonconvex sets, the normal cone--as well as the subdifferential and coderivative constructions for functions and mappings generated by it--possess comprehensive calculus rules; see \cite{rw,m06}.

Given a function $\ph\colon\R^n\to\oR:=]-\infty,\infty ]$, the (first-order) {\em subdifferential}  and {\em singular subdifferential} of $\ph$ at $\ox\in\dom\ph:=\{x\in\R^n\colon\;\ph(x)<\infty\}$ are defined geometrically via the limiting normal cone (\ref{2.2}) to the epigraph $\epi\ph:=\{(x,\al)\in\R^{n+1}\colon\;\al\ge\ph(x)\}$ by
\begin{equation}\label{2.6}
\partial\ph(\ox):=\Big\{v\in\R^n \colon\;(v,-1)\in N((\ox,\ph(\ox));{{\rm\small epi}\,\ph})\Big\},
\end{equation}
\begin{equation}\label{2.66}
\partial^\infty\ph(\ox)\colon=\Big\{v\in\R^n:\;(v,0)\in N((\ox,\ph(\ox));{{\rm\small epi}\,\ph})\Big\},
\end{equation}
respectively; the reader can find in \cite{rw,m06} useful equivalent analytic representations of (\ref{2.6}) and (\ref{2.66}). Note that, while the  subdifferential (\ref{2.6}) is a natural extension of the classical derivative and the convex subdifferential for smooth and convex functions, the singular subdifferential (\ref{2.66}) contains nonzero elements if and only if $\ph$ is {\em not} locally Lipschitzian around $\ox$ provided that it is {\em lower semicontinuous} (l.s.c.), which we {\em always assume} in what follows.

For a set-valued mapping $F\colon\R^n\rightrightarrows\R^m$ with its domain and graph
$$
\dom F:=\Big\{x\in\R^n \colon\;F(x)\ne\emp\Big\},\quad\gph F:=\Big\{(x,y)\in\R^n\times\R^m \colon\;x\in F(x)\Big\},
$$
the {\em coderivative} of $F$ at $(\ox,\oy)\in\gph F$ is defined by
\begin{equation}\label{2.8}
D^*F(\ox,\oy)(u):=\Big\{v\in\R^n \colon\;(v,-u)\in N((\ox,\oy);\gph F)\Big\}.
\end{equation}
It is an ``adjoint derivative" of set-valued mappings, which reduces to the adjoint/transposed Jacobian operator $D^*f(\ox)(u)=\{\nabla f(\ox)^*u\}$, $u\in\R^m$, for single-valued smooth mappings $F=f\colon\R^n\to\R^m$ with $\oy=f(\ox)$. Note that in the general case of nonsmooth and/or set-valued mappings the coderivative (\ref{2.8}) cannot be dual/adjoint to any tangentially generated derivative, since its values are nonconvex while the duality operation always generates convexity.

For a {\em proper} ($\dom\ph\ne\emp)$ extended-real-valued function $\ph\colon\R^n\to\oR$, recall that the {\em second-order subdifferential} of $\ph$ at $\ox\in\dom\ph$ relative to $\ov\in\partial\ph(\ox)$ is defined while following the dual ``derivative-of-derivative" approach \cite{m92} by
\begin{eqnarray}\label{2nd}
\partial^2\ph(\ox,\ov)(u)\colon=(D^*\partial\ph)(\ox,\ov)(u),\quad u\in\R^n,
\end{eqnarray}
which corresponds to the Hessian mapping $\partial^2\ph(\ox,\nabla\ph(\ox))(u)=\{\nabla^2\ph(\ox)u\}$ if $\ph$ is ${\cal C}^2$-smooth around $\ox$. Among numerous results on the second-order construction (\ref{2nd}) (see, e.g., \cite{mr,m06} and the references therein) we mention the explicit calculations of (\ref{2nd}) in terms of the initial data for the class of CPWL functions $\th\colon\R^m\to\oR$ (we use this notation for the further convenience) that can be equivalently described in one of the following ways \cite{rw}:

$\bullet$ The epigraphical set $\epi\th$ is a convex polyhedron in $\R^{m+1}$.\\

$\bullet$ There are $\alpha_i\in\R$, $l\in\N$, and $a_i\in\R^m$ for $i\in T_1\colon=\{1,\ldots,l\}$ such that $\th$ is represented by
\begin{equation}\label{eq00}
\th(z)=\left\{\begin{array}{ll}
\max\Big\{\la a_1,z\ra-\alpha_1,\ldots,\la a_l,z\ra-\alpha_l\Big\},&\mbox{if }\;z\in\dom\th,\\
\infty,&\mbox{otherwise},
\end{array}\right.
\end{equation}
where the domain set $\dom\th$ is a convex polyhedron given by
\begin{equation}\label{dom}
\dom\th=\Big\{z\in\R^m\colon\;\la d_i,z\ra\le\beta_i\;\mbox{ for all }\;i\in T_2:=\{1,\ldots,p\}\Big\}
\end{equation}
with some elements $d_i\in\R^m$, $\beta_i\in\R$, and $p\in\N$.

For simplicity we write $\th\in CPWL$ if $\th$ belongs to this class and deduce from (\ref{eq00}) that each $\th\in CPWL$ can be expressed as
\begin{equation}\label{theta}
\th(z)=\max\Big\{\la a_1,z\ra-\alpha_1,\ldots,\la a_l,z\ra-\alpha_l\Big\}+\dd(z;\dom\th),\quad z\in\R^m,
\end{equation}
via the indicator function of the domain of $\th$. Furthermore, it has been recently observed in \cite[Proposition~3.2]{ms15} that, besides (\ref{dom}), the domain of $\th$ admits the union representation $\dom\th=\bigcup^{l}_{i=1}{C_i}$ with $l$ taken from (\ref{eq00}) and with the sets $C_i$, $i\in T_1$, defined by
\begin{equation}\label{pwlr1}
C_i:=\Big\{z\in\dom\th\colon\;\la a_j,z\ra-\al_j\le\la a_i,z\ra-\al_i,\;\;\mbox{for all}\;\;j\in T_1\Big\}.
\end{equation}
Consider next the corresponding active index subsets in (\ref{pwlr1}) and (\ref{dom}) given by
\begin{equation}\label{active2}
K(\oz):=\Big\{i\in T_1\colon\;\oz\in C_i\Big\}\;\mbox{ and }\;I(\oz):=\Big\{i\in T_2\colon\;\la d_i,\oz\ra=\beta_i\Big\}
\end{equation}
and recall the formula for $\partial\th(\oz)$ at $\oz\in\dom\th$ obtained in \cite[Proposition~3.3]{ms15}:
\begin{equation}\label{fos}
\begin{array}{lll}
\partial\th(\oz)&=\conv\Big\{a_i:i\in K(\oz)\Big\}+N(\oz;\dom\th)\\&=\conv\Big\{a_i:i\in K(\oz)\Big\}+\cone\Big\{d_i\colon\;i\in I(\oz)\Big\}.
\end{array}
\end{equation}
Then for any $(\oz,\ov)\in\gph\partial\th$ we get from (\ref{fos}) that $\ov=\ov_1+\ov_1$, where
\begin{equation}\label{eq06}
\begin{array}{lll}
\disp\ov_1=\sum_{i\in K(\oz)}\bar\lm_i a_i&\mbox{ with }&\disp\sum_{i\in K(\oz)}\bar\lm_i=1,\;\bar\lm_i\ge 0,\;\mbox{ and }\\
\disp\ov_2=\sum_{i\in I(\oz)}\bar\mu_id_i&\mbox{ with }&\bar\mu_i\ge 0.
\end{array}
\end{equation}
Corresponding to (\ref{eq06}), define the index subsets of positive multipliers by
\begin{equation}\label{eq05}
J_+(\oz,\ov_1):=\Big\{i\in K(\oz)\colon\;\bar\lm_i>0\Big\},\quad J_+(\oz,\ov_2):=\Big\{i\in I(\oz)\colon\;\bar\mu_i>0\Big\}
\end{equation}
and then consider the following sets defined entirely via the parameters in (\ref{eq00}) and (\ref{dom}) along arbitrary index subsets $P_1\subset Q_1\subset T_1$ and $P_2\subset Q_2\subset T_2$:
\begin{eqnarray}\label{eq080}
\begin{array}{lll}
{\cal F}_{\tiny\{P_1,Q_1\},\{P_2,Q_2\}}:&=\span\Big\{a_i-a_j\,\colon\;i,j\in P_1\Big\}\\
&+\cone\Big\{a_i-a_j\,\colon\;(i,j)\in(Q_1\setminus P_1)\times P_1\Big\}\\
&+\span\Big\{d_i\,\colon\;i\in P_2\Big\}+\cone\Big\{d_i\,\colon\;i\in Q_2\setminus P_2\Big\},
\end{array}
\end{eqnarray}
\begin{eqnarray}\label{eq081}
\begin{array}{ll}
{\cal G}_{\tiny\{P_1,Q_1\},\{P_2,Q_2\}}:=\Big\{u\in\R^n\colon&\la a_i-a_j,u\ra=0\;\mbox{ if }\;i,j\in P_1,\\
&\la a_i-a_j,u\ra\le 0\;\mbox{ if }\;(i,j)\in(Q_1\setminus P_1)\times P_1,\\
&\la d_i,u\ra=0\;\mbox{ if }\;i\in P_2,\;\mbox{ and }\\
&\la d_i,u\ra\le 0\;\mbox{ if }\;i\in Q_2\setminus P_2\;\Big\}.
\end{array}
\end{eqnarray}
Now we are ready to formulate the precise calculation formulas for the second-order subdifferential of CPWL functions. It follows from \cite[Theorem~5.1]{ms15} that
\begin{equation}\label{2nd-val}
\begin{array}{ll}
\partial^2\th(\oz,\ov)(u)=\Big\{w\,\colon &(w,-u)\in{\cal F}_{\tiny\{P_1,Q_1\},\{P_2,Q_2\}}\times{\cal G}_{\tiny\{P_1,Q_1\},\{P_2,Q_2\}},\\
&(P_1,Q_1,P_2,Q_2)\in{\cal A}\Big\}
\end{array}
\end{equation}
for any $u\in\R^m$, where the set ${\cal A}$ of index quadruples is defined by
\begin{eqnarray}\label{eq092}
\begin{array}{ll}
{\cal A}:=\Big\{(P_1,Q_1,P_2,Q_2)\,\colon&P_1\subset Q_1\subset K,\;P_2\subset Q_2\subset I,\\
&(P_1,P_2)\in D(\oz,\ov),\;H_{\tiny\{Q_1,Q_2\}}\ne\emp\Big\}
\end{array}
\end{eqnarray}
with $K:=K(\oz)$, $I:=I(\oz)$, $H_{\tiny\{Q_1,Q_2\}}:=\{z\in\dom\th\,\colon\;K(z)=Q_1,\;I(z)=Q_2\}$,
\begin{eqnarray*}
D(\oz,\ov):=\Big\{(P_1,P_2)\subset K\times I\,\colon\;\ov\in\conv\{a_i|\;i\in P_1\}+\cone\{d_i\,\colon\;i\in P_2\}\Big\}.
\end{eqnarray*}
Furthermore, \cite[Theorem~5.2]{ms15} gives us the domain formula
\begin{equation}\label{domcod}
\begin{array}{ll}
\dom\sub^2\th(\oz,\ov)=\Big\{u\,\colon& \la a_i-a_j,u\ra=0\;\mbox{ for }\;i,j\in\Gamma(J_1),\\&\la d_t,u\ra=0\;\mbox{ for }\;t\in\Gamma(J_2)\Big\},
\end{array}
\end{equation}
where the index sets $\Gamma(J_1)$ and $\Gamma(J_2)$ are defined by
\begin{equation}\label{feature}
\begin{array}{ll}
\Gamma(J_1):=\Big\{i\in K\,\colon\;\la a_i-a_j,u\ra=0\;\mbox{ for all }\;j\in J_1\;\mbox{ and }\;u\in{\cal G}_{\tiny\{J_1,K\},\{J_2,I\}}\;\Big\},\\
\Gamma(J_2):=\Big\{t\in I\,\colon\;\la d_t,u\ra=0\;\mbox{ for all }\;u\in{\cal G}_{\tiny\{J_1,K\},\{J_2,I\}}\;\Big\}
\end{array}
\end{equation}
with the notation $J_1:=J_+(\oz,\ov_1)$ and $J_2:=J_+(\oz,\ov_2)$ from (\ref{eq06}), (\ref{eq05}).

In the subsequent sections of the paper, we will consider compositions $\th\circ\Phi$ of CPWL outer functions $\th\colon\R^m\to\oR$ and inner mappings $\Phi\colon\R^n\times\R^d\to\R^m$ that are ${\cal C}^2$-smooth around some $(\ox,\ow)$ with $\oz:=\Phi(\ox,\ow)\in\dom\th$ under the first-order qualification condition
\begin{equation}\label{2.9}
\sub^{\infty}\th(\oz)\cap\ker\nabla_x\Phi(\ox,\ow)^*=\{0\}.
\end{equation}
Such compositions form an important subclass of functions known as {\em fully amenable} in $x$ at $\ox$ with compatible parametrization by $w$ at $\ow$ (we will drop in what follows the latter parametrization expression for brevity), which are defined in this way with using more general convex {\em piecewise linear-quadratic} outer functions $\th$; see \cite{rw} for more details.\vspace*{-0.2in}

\section{Reducibility, Nondegeneracy and Second-Order Qualification}\label{nondeg}

The main goal of this section is to establish relationships between the second-order qualification condition introduced in \cite{rw} in order to derive the exact second-order chain rule for fully amenable compositions with CPWL outer functions and the partial nondegeneracy condition of a completely different nature that was employed in \cite{mnn} to get the same second-order chain rule. In this way we obtain below some auxiliary results of their independent interest.

Considering first a fully amenable composition $\psi=\th\circ\Phi$ as defined at the end of Section~\ref{prel}, recall that the {\em second-order qualification condition} (SOCQ) holds for $\psi$ in $x$ at $(\ox,\ow)$ if
\begin{equation}\label{2qcf}
\partial^2\th(\oz,v)(0)\cap\ker\nabla_x\Phi(\ox,\ow)^*=\{0\}\;\;\mbox{for all}\;v\in M(\ox,\ow,\oq),
\end{equation}
where $\oq\in\partial_x\psi(\ox,\ow)$ is a fixed partial subgradient of $\psi$ in $x$ at $(\ox,\ow)$ and
\begin{eqnarray}\label{M1}
M(\ox,\ow,\oq):=\Big\{v\in\R^m\,\colon\;v\in\partial\th(\oz)\;\mbox{ with }\;\nabla_x\Phi(\ox,\ow)^*v=\oq\Big\}.
\end{eqnarray}
Note that the imposed qualification condition (\ref{2.9}) ensures, by using the well-known first-order subdifferential chain rule \cite{m06,rw}, that $M(\ox,\ow,\oq)\ne\emp$.

For a given $\th\colon\R^n\to\oR$, denote by $S(z)$ a subspace of $\R^m$ parallel to the {\em affine hull} $\aff\partial\th(z)$ of the subdifferential $\partial\th(z)$, $z\in\R^m$. The next theorem provides a precise calculation of the second-order subdifferential for CPWL functions at the {\em origin} $0\in\R^m$ entirely via the initial data in (\ref{eq00}) and (\ref{dom}), relates it to the subspace $S(\oz)$ defined above, and gives an effective representation of SOQC in (\ref{2qcf}) convenient for our further analysis and applications.

\begin{theorem} {\bf(second-order subdifferential of CPWL functions at the origin and SOQC representation).}\label{sospwl}
Let $\psi=\th\circ\Phi$ be a fully amenable composition of $\th\colon\R^m\to\oR$ and $\Phi\colon\R^n\times\R^d\to\R^m$ with $\th\in CPWL$, and let $S(\oz)$ be a subspace of $\R^m$ parallel to the affine hull $\aff\partial\th(\oz)$ with $\oz=\Phi(\ox,\ow)$. Then the following assertions hold:

{\bf(i)} For all $\ov\in\partial\th(\oz)$ we have the representation
\begin{eqnarray}\label{2nd0}
\sub^2\th(\oz,\ov)(0)=\span\Big\{a_i-a_j\,\colon\;i,j\in K(\oz)\Big\}+\span\Big\{d_i\,\colon\;i\in I(\oz)\Big\}
\end{eqnarray}
via the data in {\rm(\ref{eq00})}, {\rm(\ref{dom})} with the active index sets $K(\oz)$, $I(\oz)$ defined in {\rm(\ref{active2})}.

{\bf(ii)} Furthermore, we have $\sub^2\th(\oz,\ov)(0)=S(\oz)$ independently of $\ov\in\partial\th(\oz)$.

{\bf(iii)} The SOQC property {\rm(\ref{2qcf})} can be equivalently written as
\begin{equation}\label{2qcf2}
S(\oz)\cap\ker\nabla_x\Phi(\ox,\ow)^*=\{0\}.
\end{equation}
independently of $\ov\in\partial\th(\oz)$ and $\oq\in\partial_x\psi(\ox,\ow)$ in {\rm(\ref{M1})}.
\end{theorem}
{\it Proof} To verify first the inclusion ``$\subset$" in (\ref{2nd0}), pick $y\in\sub^2\th(\oz,\ov)(0)$ and find by (\ref{2nd-val}) an index quadruple $(P_1,Q_1,P_2,Q_2)\in{\cal A}$ from (\ref{eq092}) such that
$$
(y,0)\in{\cal F}_{\tiny\{P_1,Q_1\},\{P_2,Q_2\}}\times{\cal G}_{\tiny\{P_1,Q_1\},\{P_2,Q_2\}}.
$$
We immediately deduce from representation (\ref{eq080}) that
$$
{\cal F}_{\tiny\{P_1,Q_1\},\{P_2,Q_2\}}\subset\span\Big\{a_i-a_j\,\colon\;i,j\in K(\oz)\Big\}+\span\Big\{d_i\,\colon\;i\in I(\oz)\Big\},
$$
which justifies the inclusion ``$\subset$" in (\ref{2nd0}). To derive further the opposite inclusion ``$\supset$" therein, take any vector $y$ with
$$y\in\span\Big \{a_i-a_j\,\colon\;i,j\in K(\oz)\Big\}+\span\Big\{d_i\,\colon\;i\in I(\oz)\Big\}$$
and  then put $P_1=Q_1:=K(\oz)$ and $P_2=Q_2:=I(\oz)$. Since $\oz\in H_{\{Q_1,Q_2\}}$ in (\ref{2nd0}), it follows that $(P_1,Q_1,P_2,Q_2)\in{\cal A}$. Employing again the second-order formula (\ref{2nd-val}) tells us that $(y,0)\in N((\oz,\ov);\gph\partial\th)$ and hence yields $y\in\sub^2\th(\oz,\ov)(0)$, which thus verifies assertion (i).

To prove assertion (ii), observe that $S(\oz)=\aff\partial\th(\oz)-a_t$ for some $t\in K(\oz)$. Picking $y\in S(\oz)$ gives us $y+a_t=\sum_{i=1}^{s}\alpha_ic_i$ for some vectors $c_i\in\sub\th(\oz)$ and some number $s>0$ with $\sum_{i=1}^{s}\alpha_i=1$. It follows from (\ref{fos}) that $c_i=c_{1i}+c_{2i}$ with $c_{1i}\in \co\{a_r\,\colon\;r\in K(\oz)\}$ and $c_{2i}\in N(\oz;\dom\th)$ for $i=1,\ldots,s$. Therefore
\begin{equation}\label{eq0450}
y=\sum_{i=1}^{s}\alpha_ic_i-a_t=\sum_{i=1}^s\alpha_i(c_{1i}-a_t)+\sum_{i=1}^s\alpha_i c_{2i}.
\end{equation}
It is clear that $c_{1i}-a_t\in\span\{a_i-a_j\,\colon\;i,j\in K(\oz)\}$, and thus we have
$$
c_{i}-a_t\in\span\Big\{a_i-a_j\,\colon\;i,j\in K(\oz)\Big\}+\span\Big\{d_i\,\colon\;i\in I(\oz)\Big\}.
$$
Using this together with (\ref{eq0450}), (\ref{2nd0}) justifies $S(\oz)\subset\sub^2\th(\oz,\ov)(0)$ in (ii).

To verify the opposite inclusion, take $y\in\sub^2\th(\oz,\ov)(0)$ and get by (\ref{2nd0}) that
$$
y=\sum_{(i,j)\in A_1\times A_2}\alpha_{i,j}(a_i-a_j)+\sum_{t\in A_3}\beta_t d_t
$$
with some index subsets $A_1,A_2\subset K(\oz)$ and $A_3\subset I(\oz)$. Select now the subsets $B_1,B_2\subset A_1\times A_2$ and $B_3,B_4\subset A_3$ so that
\begin{equation}\label{eq0451}
\begin{array}{lll}
A_1\times A_2=B_1\cup B_2&\mbox{and}&A_3=B_3\cup B_4,\\
\al_{i,j}\ge 0\;\;\mbox{whenever}\;\;(i,j)\in B_1&\mbox{and}&\al_{i,j}<0\;\;\mbox{whenever}\;\;(i,j)\in B_2,\\
\beta_{t}\ge 0\;\;\;\;\mbox{whenever }\;\;t\in B_3&\mbox{and}&\beta_{t}<0\;\;\;\;\mbox{whenever }\;\;t\in B_4.
\end{array}
\end{equation}
In this way we represent the given vector $y$ as $y=y'-b$ with
$$
\begin{array}{lll}
\disp y':=\sum_{(i,j)\in B_1}\alpha_{i,j}a_i-\sum_{(i,j)\in B_2}\alpha_{i,j}a_j+\sum_{t\in B_3}\beta_t d_t\;\mbox{ and }\\
\disp b:=\sum_{(i,j)\in B_1}\alpha_{i,j}a_j-\sum_{(i,j)\in B_2}\alpha_{i,j}a_i+\sum_{t\in B_4}(-\beta_t)d_t.
\end{array}
$$
Denoting $\al:=\sum_{(i,j)\in B_1}\alpha_{i,j}-\sum_{(i,j)\in B_2}\alpha_{i,j}$, deduce from (\ref{eq0451}) that $\al\ge 0$. For $\al>0$  we get the inclusions
\begin{eqnarray}\label{Sb}
\frac{1}{\al}y'\in\aff\partial\th(\oz)\;\mbox{ and }\;\frac{1}{\al}b\in\aff\partial\th(\oz).
\end{eqnarray}
It follows from the construction of $S(\oz)$ and the second inclusion in (\ref{Sb}) that  $S(\oz)=\aff\partial\th(\oz)-\frac{1}{\al}b$, and so the first one in (\ref{Sb}) yields $\frac{1}{\al}y\in S(\oz)$. This shows that $y\in S(\oz)$ since $S(\oz)$ is a subspace, and thus we get $\sub^2\th(\oz,\ov)(0)\subset S(\oz)$ in the case of $\al>0$. Considering now the remaining case of $\al=0$ gives us the expression $y=\sum_{t\in A_3}\beta_t d_t$, which implies by (\ref{eq0451}) that
$$
y=a_t+\sum_{t\in B_3}\beta_t d_t-\Big(a_t+\sum_{t\in B_4}(-\beta_t)d_t\Big)\in S(\oz)\;\mbox{ with some }\;t\in K(\oz)
$$
due to $a_t+\sum_{t\in B_4}(-\beta_t)d_t\in\aff\partial\th(\oz)$ and hence verifies (ii). This immediately yields (\ref{2qcf2}) in (iii) by comparison it with the SOCQ definition in (\ref{2qcf}). $\b$

Note that while the precise calculation of $\sub^2\th(\oz,\ov)(0)$ in Theorem~\ref{sospwl}(i) is new, assertion (ii) therein follows from the proof of Theorem~4.3 in \cite{mr} by using the representation of $\sub^2\th(\oz,\ov)(0)$ for piecewise linear-quadratic functions $\th\colon\R^m\to\oR$ established in \cite[Theorem~4.1]{mr}. The proof of the latter result in \cite{mr} is based on the tangential approach from \cite{rw} being significantly more involved in comparison with the one given above.

It is also worth mentioning as a by-product of the above calculations that the validity of SOQC for fully amenable compositions with CPWL outer functions yields the fulfillment of the first-order qualification condition (\ref{2.9}) in the definition of such compositions. To see this, recall that $\partial^\infty\th(\oz)=N(\oz;\dom\th)$ for convex functions and thus get $\partial^\infty\th(\oz)\subset\sub^2\th(\oz,\ov)(0)$ whenever $\ov\in\partial\th(\oz)$ by comparing (\ref{2nd0}) with that of $N(\oz;\dom\th)=\cone\{d_i\,\colon\;i\in I(\oz)\}$.

Next we consider the concept of {\em nondegeneracy}. It was first initiated for {\em sets} in \cite{rob84} as a polyhedral counterpart of the classical LICQ in nonlinear programming. Note that even for MPPCs this nondegeneracy condition may be strictly weaker than LICQ; see \cite{mrs} for equivalent descriptions for MPPCs and particularly Example~6.7 therein. Nondegeneracy and associated {\em reducibility} notions for general sets were comprehensively studied in \cite{bs} based on the previous paper of these authors. For the case of extended-real-valued {\em functions} the notion of ${\cal C}^2$-reducibility and the corresponding notion of partial nondegeneracity was formulated in \cite{mnn} in order to derive the aforementioned second-order subdifferential chain rule; see below.

Following this pattern, we say that a function $\th\colon\R^m\to\oR$ is {\em ${\cal C}^2$-reducible} (resp.\ {\em ${\cal C}^\infty$-reducible}) to a function $\vt\colon\R^s\to\oR$ at $\oz$ with $s\le m$ if there exists a ${\cal C}^2$-smooth (resp.\ ${\cal C}^\infty$-smooth) mapping $h\colon\R^m\to\R^s$ with the surjective derivative $\nabla h(\oz)$ such that $\th(z)=(\vartheta\circ h)(z)$ for all $z$ around $\oz$.

Our next result shows that any function $\th\in CPWL$ on $\R^m$ is ${\cal C}^\infty$-reducible to some $\vt\in CPWL$ on $\R^s$ by using a {\em linear} surjective operator $h\colon\R^m\to\R^s$. From now on we assume that $0\in\aff\partial\th(\oz)$ at $\oz\in\dom\th$, which implies that  $S(\oz)=\aff\partial\th(\oz)$. In fact this assumption does {\em not restrict the generality} in dealing with the second-order subdifferential. Indeed, we have $S(\oz)=\aff\partial\theta(\oz)-b_{\oz}$ for some $b_{\oz}\in\aff\partial\theta(\oz)$. Defining then
$\bar{\theta}(z):=\theta(z)-\la b_{\bar{z}},z\ra$ shows that $0\in\aff\partial\bar\th(z)$ and $\partial^2\th(\oz,\oy)=\partial^2\bar\th(\oz,\oy-b_{\oz})$ for any $\ov\in\partial\th(\oz)$.

\begin{lemma}{\bf(${\cal C}^\infty$-reducibility of  piecewise linear functions).}\label{fred} Let the function $\th\colon\R^m\to\oR$ be CPWL, let $\oz\in\dom\th$, and let $s:=\dim S(\oz)\le m$. Then $\th$ is ${\cal C}^\infty$-reducible at $\oz$ to a CPWL function $\vt\colon\R^s\to\oR$ via a linear operator $h(z):=Bz$ generated by some $s\times m$ matrix $B$.
\end{lemma}
{\it Proof} It follows from \cite[Proposition~3.3(i)]{ms15} that $\partial\th(z)\subset\partial\th(\oz)$ for all $z\in O$ in some neighborhood of $\oz$.
Denote by $A$ the matrix of a linear isometry from $\mathbb{R}^m$ into $\mathbb{R}^s\times\mathbb{R}^{m-s}$ under which $A^*(S(\oz))=\mathbb{R}^s\times\{0\}$. Define $\xi\colon\R^m\to\oR$ by
\begin{eqnarray}\label{xi}
\xi(y):=\th(Ay)\;\mbox{ for all }\;y\in\R^m
\end{eqnarray}
and get by \cite[Proposition~3.55(b)]{rw} that $\xi$ is proper, convex, and piecewise linear on $\R^m$. Applying the chain rule of convex analysis to (\ref{xi}) gives us
\begin{equation}\label{5.5}
\partial\xi(y)=A^*\partial\theta(z)\;\mbox{ with }\;Ay=z.
\end{equation}
Denote $U:=A^{-1}(O)$ and deduce from the classical open mapping theorem that $U$ is a neighborhood of $\oy:=A^{-1}\oz$. Suppose that $\alpha>0$ is so small that $\B_{\alpha}(\oy)\subset U$ for the ball centered at $\oy$ with radius $\al$. Letting  $O':=A(\int\B_{\alpha}(\oy))$, we deduce from the open mapping theorem that $O'$ is a neighborhood of $\oz$. Then $S(z)=\aff\partial\theta(z)+b_z$ with some $b_z\in\mathbb{R}^m$ for each $z\in O$, and furthermore $b_{\oz}=0$ as discussed before the formulation of the lemma. This tells us that
\begin{equation}\label{s4a}
v=(v_1,\ldots,v_m)\in\partial\xi(y)=A^*\partial\theta(z)\subset A^*\partial\theta(\bar{z})\subset A^*(S(\oz))-A^*b_{\bar{z}}\subset\mathbb{R}^s\times\{0\}
\end{equation}
for all $y\in U$, which implies that the last $m-s$ elements of any $v\in\partial\xi(y)$ are zeros whenever $y\in U$. Construct now the desired $s\times m$ matrix $B$ claimed in the lemma from the $m\times m$ matrix $A^{-1}$ by deleting the last $m-s$ rows of the latter. We define the corresponding function $\vt\colon\R^s\to\oR$ by using $\xi$ in (\ref{xi}) as follows: take $y=(y_s,y_{m-s})=(x,y_{m-s})\in\R^s\times\mathbb{R}^{m-s}$ and put
\begin{eqnarray}\label{vt}
\vt(x):=\xi(x,\bar y_{m-s})=\xi(y_s,\bar y_{m-s})\;\mbox{ for all }\;x\in\R^s,
\end{eqnarray}
where $\bar y_{m-s}$ is the last $m-s$ elements of the vector $\bar y=A^{-1}\oz$. Since $\xi$ is proper, so is the function $\vt$ in (\ref{vt}). It is easy to see that $\vt$ is piecewise linear and the convexity of $\xi$ implies the convexity of $\vartheta$. To justify the statement of the lemma, it remains to verify the representation
\begin{eqnarray}\label{reduc}
\th(z)=(\vt\circ B)(z)\;\mbox{ for all }\;z\in O'.
\end{eqnarray}
Let us do it by observing first that $y\in\int\B_{\alpha}(\oy)$ whenever $y=A^{-1}z$ generated by $z\in O'$. It follows from (\ref{xi}), (\ref{vt}), and the definition of the matrix  $B$ that $(\vartheta\circ B)(z)=\xi(y_s,\bar y_{m-s})$ in the notation above, where $(y_s,\bar y_{n-s})\in\int\B_{\alpha}(\oy)$. Thus (\ref{reduc}) would be implied by the relationship
\begin{eqnarray}\label{reduc1}
\xi(y_s,\bar y_{m-s})=\xi(y_s,y_{m-s})\;\mbox{ for any }\;y=(y_s,y_{m-s})=A^{-1}z,\quad z\in O'.
\end{eqnarray}
Since (\ref{reduc1}) is trivial when both values $\xi(y_s,y_{m-s})$ and $\xi(y_s,\oy_{m-s})$ are infinite, suppose without loss of generality that $\xi(y_s, y_{n-s})$ is a real number. The polyhedrality of $\epi\xi$ ensures that the function $\xi$ is l.s.c., and hence we can apply to it the approximate mean value inequality from \cite[Corollary~3.50]{m06}. This allows us to find $c\in\R^m$ on the segment connecting $(y_s,y_{m-s})$ and $(y_s,\bar y_{m-s})$ as well as a sequence $v_k\in\partial\xi(u_k)$ with $u_k\to c$ and $\xi(u_k)\to\xi(c)$ so that
\begin{equation}\label{5.4}
\xi(y_s,\bar y_{m-s})-\xi(y_s,y_{m-s})\le\liminf_{k\to\infty}\la v_k,(0_s,\bar y_{m-s}-y_{m-s})\ra.
\end{equation}
It follows from (\ref{s4a}) that $u_k\in\int\B_{\alpha}(\oy)\subset U$ and so $\la v_k,(0_s,\bar y_{m-s}-y_{m-s})\ra=0$ for all $k\in\N$ sufficiently large. In the same way we get the opposite inequality
\begin{eqnarray*}
\xi(y_s,y_{m-s})-\xi(y_s,\bar y_{m-s})\le 0
\end{eqnarray*}
and combining it with (\ref{5.4}) arrive at (\ref{reduc1}), which completes the proof. $\b$

Now we are ready to formulate, following \cite{mnn}, the notion of {\em nondegeneracy} of one mapping relative to another one used for deriving the second-order chain rule. Observe that, although this notion is formulated for two arbitrary mappings, its application to second-order analysis mainly concerns amenable compositions $\th\circ\Phi$ of $\th\colon\R^m\to\oR$ and $\Phi\colon\R^n\times\R^d\to\R^m$  while defining nondegenerate points of $\Phi\colon\R^n\times\R^d$ relative to the mapping $h\colon\R^m\to\R^s$ that furnishes the appropriate {\em reducibility} of the outer function $\th$. Thus in our case of $\th\in CPWL$ we deal with {\em linear} mapping $h(z)=Bz$ that appears in the ${\cal C}^\infty$-reducibility assertion of Lemma~\ref{fred}.

Having this in mind, it is said that $(\ox,\ow)\in\R^n\times\R^d$ is a {\em partial nondegenerate point} of $\Phi\colon\R^n\times\R^d\to\R^m$ in $x$ relative to $h\colon\R^m\to\R^s$ if
\begin{equation}\label{fnond}
\nabla_x\Phi(\ox,\ow)\R^n+\ker\nabla h(\oz)=\R^m\;\mbox{ with }\;\oz=\Phi(\ox,\ow)
\end{equation}
under the corresponding differentiability assumptions on $\Phi$ and $h$. The next theorem based on the previous results of this section reveals that, in the case of fully amenable compositions with CPWL outer functions, the SOQC property (\ref{2qcf2}) of $\th\circ\Phi$ is {\em equivalent} to the nondegeneracy condition (\ref{fnond}) {\em provided} that the mapping $h\colon\R^m\to\R^s$ with $s=\dim S(\oz)$ therein is the linear transformation $h(z)=Bz$ constructed in Lemma~\ref{fred} to realize the ${\cal C}^\infty$-reducibility of $\th$.

It is worth mentioning that this line of equivalency between the corresponding SOQC and nondegeneracy properties is a continuation of the results previously established in \cite{mrs} in connection with mathematical programs with polyhedral constraints and in \cite{mos} in connection with second-order cone programs (SOCPs), where (in both cases) the nondegeneracy condition of a mapping relative to the underlying set (polyhedron and second-order cone, respectively) was understood in the sense of \cite{bs} via the tangent cone to this set. The crucial difference of our case in this paper is that we implement the general nondegeneracy/reducibility notion \cite{mnn} relative to a mapping and emphasize the {\em linearity} of this mapping in the CPWL setting under consideration.

\begin{theorem}{\bf(relationship between SOQC and nondegeneracy for fully amenable compositions with CPWL outer functions).}\label{fred2} Consider a fully amenable composition $\psi=\th\circ\Phi$, which is finite at $(\ox,\ow)\in\R^n\times\R^d$. Let $\th\colon\R^m\to\oR$ be CPWL and let $B$ be an $s\times m$ matrix constructed in Lemma~{\rm\ref{fred}}. Then the SOQC property {\rm (\ref{2qcf2})} holds at $(\ox,\ow)$ if and only if this point is partially nondegenerate {\rm(\ref{fnond})} for $\Phi$ relative to $h(z)=Bz$ with $s=\dim S(\oz)$.
\end{theorem}
{\it Proof} Since $0\in\aff\partial\th(\oz)$ as discussed before the formulation of Lemma~\ref{fred}, we have $S(\oz)=\aff\partial\th(\oz)$. This lemma gives us a CPWL function $\vt\colon\R^s\to\oR$ and a mapping $h(z)=Bz$ from $\R^m$ to $\R^s$ such that $\th(z)=(\vartheta\circ h)(z)$ for all $z\in\R^m$ sufficiently close to $\oz$. Assuming that SOQC holds at $(\ox,\ow)$ and taking the orthogonal complements of both sides in (\ref{2qcf2}), we arrive at
\begin{equation}\label{5.10}
\nabla_x\Phi(\ox,\ow)\R^n+S(\oz)^\bot=\R^m.
\end{equation}
To deduce from (\ref{5.10}) the partial nondegeneracy condition (\ref{fnond}) with $h(z)=Bz$, it suffices to show that $\ker\nabla h(\oz)= S(\oz)^\bot$, which reads as
$\ker B=S(\oz)^\bot$. Indeed, picking $u\in\ker B$ and taking into account that $A^*(S(\oz))=\R^s\times\{0\}$ in the proof of the lemma yield
$$
0=\la A^{-1}u,A^*p\ra=\la u,(A^{-1})^*A^*p\ra=\la u,p\ra\;\;\mbox{for any}\;\;p\in S(\oz),
$$
which tells us that $u\in S(\oz)^\bot$, and so $\ker B\subset S(\oz)^\bot$. The opposite inclusion $S(\oz)^\bot\subset\ker B$ can be checked similarly, which shows therefore that SOQC$\Longrightarrow$partial nondegeneracy. The same arguments allow us to verify via (\ref{5.10}) that  partial nondegeneracy$\Longrightarrow$SOQC, which completes the proof. $\b$

The final result of this section presents the second-order chain rule for the {\em partial} second-order subdifferential (denoted below as $D^*\partial_x\psi$) of fully amenable compositions $\psi=\th\circ\Phi$ with CPWL outer functions. This result was first obtained in \cite[Theorem~4.3]{mr} for nonparametric compositions and then in \cite[Theorem~4.1]{mrs} in the general parametric case. Both proofs in \cite{mr,mrs} are involved being based on the difficult Theorem~4.1 from \cite{mr}. The new proof given below is much simpler based on the equivalency result of Theorem~\ref{fred2} and the second-order chain rule obtained in \cite[Theorem~3.6]{mnn} under nondegeneracy condition in the Banach setting.

\begin{corollary}{\bf (second-order chain rule for parametric compositions with CPWL outer functions).}\label{soc} Let $\psi=\th\circ\Phi$ be a fully amenable composition with a CPWL function $\th\colon\R^m\to\oR$ and an inner mapping $\Phi\colon\R^n\times\R^d\to\R^m$ that is
${\cal C}^2$-smooth around $(\ox,\ow)$. Then the validity of SOQC in {\rm (\ref{2qcf2})} ensures that for any $\oq\in\partial_x\psi(\ox,\oy)$ the set $M(\ox,\ow,\oq)$ from {\rm(\ref{M1})} is a singleton $\{\ov\}$ and we have the following second-order chain rule whenever $u\in\R^n$:
\begin{equation}\label{2cr}
\begin{array}{ll}
(D^*\partial_x\psi)(\ox,\ow,\oq)(u)&=\Big(\nabla^2_{xx}\la\ov,\Phi\ra(\ox,\ow)u,\nabla^2_{xw}\la\ov,\Phi\ra(\ox,\ow)u\Big)\\
&+\Big(\nabla_x\Phi(\ox,\ow),\nabla_w\Phi(\ox,\ow)\Big)^*\partial^2\th(\oz,\oq)(\nabla_x\Phi(\ox,\ow)u).
\end{array}
\end{equation}
\end{corollary}
{\it Proof.} For any $\ov_1,\ov_2\in M(\ox,\ow,\oq)$ we have $\ov_1-\ov_2\in\ker\nabla_x\Phi(\ox,\ow)^*$. Because $\ov_1,\ov_2\in\sub\th(\oz)$, it follows from (\ref{fos}) and Theorem~\ref{sospwl}(i) that $\ov_1,\ov_2\in S(\oz)$. Applying now SOCQ (\ref{2qcf2}) gives us $\ov_1=\ov_2$, and so $M(\ox,\ow,\oq)=\{\ov\}$. Then we get from Lemma~\ref{fred} that $\th$ is ${\cal C}^\infty$-reducible by the linear mapping $h(z)=Bz$, and hence $(\ox,\ow)$ is a partial nondegenerate point (\ref{fnond}) of $\Phi$ relative to this mapping $h\colon\R^m\to\R^s$ with $s=\dim S(\oz)$. To arrive finally at the chain rule (\ref{2cr}), it remains to apply \cite[Theorem~3.6]{mnn} and thus complete the proof. $\b$

Note that Corollary~\ref{soc} clarifying the meaning of \cite[Theorem~4.3]{mr} and \cite[Theorem~4.3]{mrs} can be viewed as a realization of the second-order chain rule from \cite[Theorem~3.6]{mnn} in the case of CPWL outer functions under the fulfillment of SOQC, which corresponds to a {\em linear reduction mapping} $h\colon\R^n\to\R^s$ in the nondegeneracy condition (\ref{fnond}). The result of the latter theorem justifies the validity of (\ref{2cr}) under (\ref{fnond}) when $h$ is merely a {\em ${\cal C}^2$-smooth} mapping that furnishes the required reducibility of $\th$. In what follows we refer to the second-order chain rule (\ref{2cr}) valid under the nondegeneracy condition (\ref{fnond}) with some ${\cal C}^2$-smooth reduction mapping $b$.\vspace*{-0.2in}

\section{Full Stability in Composite Optimization}\label{comp}

In this section, we proceed with applications of second-order generalized differentiation to problems of {\em composite optimization} given in the form:
\begin{equation}\label{fcp1}
\mbox{minimize }\ph_0(x)+\theta(\Phi(x))\;\mbox{s.t.}\;x\in\R^n\;\mbox{with}\;\Phi(x):=\big(\ph_1(x),\ldots,\ph_m(x)\big),
\end{equation}
where $\theta\colon\R^m\to\Bar\R$ is a CPWL extended-real-valued function, and where all $\ph_i\colon\R^n\to\R$, $i=0,\ldots,m$, are ${\cal C}^2$-smooth around the reference optimal solution. This class of problems encompasses conventional  problems of nonlinear programming (NLPs), MPPCs mentioned in Section~3 as well as constrained and unconstrained minimax problems. It also includes the following major subclass of {\em extended nonlinear programs} (ENLPs) introduced in \cite{r}:
\begin{equation}\label{enlp}
\mbox{minimize }\;\ph_0(x)+(\th\circ\Phi)(x)\;\mbox{ with }\;\th(z):=\disp\sup_{p\in P}\langle p,z\rangle,\quad x\in\R^n,
\end{equation}
where $P$ is a convex polyhedron, and so $\th$ in (\ref{enlp}) is piecewise linear; see \cite{ms15}.

Consider now the two-parametric version of (\ref{fcp1}) constructed by
\begin{equation}\label{fcp2}
{\cal P}(w,v):\quad\mbox{minimize }\;\ph_0(x,w)+\theta(\Phi(x,w))-\la v,x\ra\;\mbox{ subject to }\;x\in\R^n,
\end{equation}
where the perturbed functions $\ph_0(x,w)$ and $\Phi(x,w)=(\ph_1(x,w),\ldots,\ph_m(x,w))$ are ${\cal C}^2$-smooth with respect to both variables. Denote
\begin{equation}\label{phi}
\ph(x,w):=\ph_0(x,w)+\th(\Phi(x,w))\;\mbox{ for }\;(x,w)\in\R^n\times\R^d
\end{equation}
and then fix $\gg>0$ and $(\ox,\ow,\ov)$ with $\Phi(\ox,\ow)\in\dom\th$ and $\ov\in\partial_x\ph(\ox,\ow)$. Define the parameter-depended  optimal value function for (\ref{fcp2}) by
\begin{eqnarray*}
m_\gg(w,v):=\inf_{\|x-\ox\|\le\gg}\Big\{\ph(x,w)-\la v,x\ra\Big\}
\end{eqnarray*}
and the parameterized set of optimal solutions to (\ref{fcp1}) by
\begin{eqnarray}\label{M}
M_\gg(w,v):=\mbox{argmin}_{\|x-\ox\|\le\gg}\Big\{\ph(x,w)-\la v,x\ra\Big\}
\end{eqnarray}
with the convention that argmin:=$\emp$ when the expression under minimization is $\infty$. Following \cite{lpr}, we say that a point $\ox$ is a {\em fully stable} locally optimal solution to problem ${\cal P}(\ow,\ov)$ in (\ref{fcp2})  if there exist a number $\gg>0$ and neighborhoods $W$ of $\ow$ and $V$ of $\ov$ such that the mapping $(w,v)\mapsto M_\gg(w,v)$ is single-valued and Lipschitz continuous with $M_\gg(\ow,\ov)=\{\ox\}$ and the function $(w,v)\mapsto m_\gg(w,v)$ is likewise Lipschitz continuous on $W\times V$.

Full stability of local minimizers was initiated and characterized in \cite{lpr} in the extended-real-valued format of unconstrained optimization in finite-dimensional spaces. Recent second-order characterizations of full stability have been obtained under various constraint qualifications for NLPs, ENLPs, and MPPCs \cite{mrs}, SOCPs \cite{mos} mentioned in Section~3, problems of polyhedric programming in Hilbert spaces with applications to optimal control of semilinear PDEs \cite{mn3}, general conic programs with applications to semidefinite programming \cite{mnr}, and unconstrained minimax problems \cite{ms14}. Furthermore, this notion has been extended to parametric variational systems \cite{mn14} with second-order characterizations and applications to variational inequalities in Hilbert spaces and variational conditions in finite dimensions.

In this section we establish new second-order characterizations of full stability for local optimal solutions to problems of composite optimization (\ref{fcp2}) with CPWL outer functions therein. In particular, the results established below cover those in \cite{mrs,ms14} while being independent from characterizations obtained in \cite{mos,mn3,mnr,mn14} for optimization and variational problems that cannot be represented in the composite form (\ref{fcp1}) with a CPWL outer function $\th$.

To proceed, denote $\oz:=\Phi(\ox,\ow)\in\dom\th$ and recall from Lemma~\ref{fred} that $\th$ is reducible at $\oz$ to some CPWL function $\vt\colon\R^s\to\oR$ with $s=\dim S(\oz)$ by using a linear mapping $h(z)=Bz$ with the $s\times m$ matrix $B$ constructed in that lemma. Thus we have $\th(z)=(\vartheta\circ B)(z)$ for all $z$ near $\oz$ generating the mapping $\Psi(x,w):=(B\circ\Phi)(x,w)$. This tells us that the problem
${\cal P}(w,v)$ from (\ref{fcp2}) is locally equivalent around  $(\ox,\ow)$ to the following {\em reduced} problem:
\begin{equation}\label{fcp3}
{\cal P}_r(w,v):\quad\mbox{minimize }\;\ph_0(x,w)+\vartheta(\Psi(x,w))-\la v,x\ra\;\mbox{ subject to }\;x\in\R^n.
\end{equation}
We will see below that the reduced problem (\ref{fcp3}) is instrumental in deriving the explicit second-order characterization of full stability of local minimizers in composite optimization obtained in this section as well as other important results established later on. The main assumption we need in what follows is the following {\em nondegeneracy condition} discussed in Section~3:
\begin{itemize}
\item [{\bf ND:}] A pair $(\ox,\ow)$ is a  partial nondegenerate point (\ref{fnond}) of $\Phi$ from (\ref{fcp2}) in $x$ relative to the linear mapping
$h(z)=Bz$, where $B$ is the $s\times m$ matrix constructed in the proof of Lemma~\ref{fred} with $s=\dim S(\oz)$, $\oz=\Phi(\ox,\ow)$.
\end{itemize}
We know from Theorem~\ref{fred2} that condition ND is equivalent to the SOCQ property (\ref{2qcf2}) in the framework of the composite optimization problem (\ref{fcp2}).

The next proposition is a composite optimization counterpart of \cite[Proposition~3.5]{ms14} obtained therein for constrained optimization problems with $\th=\dd_\Th$, the indicator function of a ${\cal C}^2$-reducible closed and convex set $\Th$.

\begin{proposition}{\bf(full stability and nondegeneracy in the original and reduced problems).}\label{red1} Let $\ox$ be a feasible solution to ${\cal P}(\ow,\ov)$ from {\rm(\ref{fcp2})} for the parameter pair $(\ow,\ov)\in\R^d\times\R^n$. Then $\ox$ is a fully stable locally optimal solution to ${\cal P}(\ow,\ov)$ if and only if it is a fully stable locally optimal solution to ${{\cal P}_r}(\ow,\ov)$. Furthermore, the validity  of {\rm ND} for $(\ox,\ow)$ implies the surjectivity $($full rank$)$ of the partial Jacobian matrix $\nabla_x\Psi(\ox,\ow)$, where $\Psi=B\circ\Phi$.
\end{proposition}
{\it Proof} The claimed equivalence follows directly from the above observation that problems ${\cal P}(w,v)$ and ${\cal P}_r(w,v)$ are locally the same. Let us verify the part of the proposition concerning nondegeneracy. Supposing that ND holds gives us
$$
\nabla_x\Phi(\ox,\ow)\R^n+\ker\nabla h(\oz)=\nabla_x\Phi(\ox,\ow)\R^n+\ker B=\R^m.
$$
It yields by applying the classical chain rule that
$$
\nabla_x\Psi(\ox,\ow)\R^n=B\nabla_x\Phi(\ox,\ow)\R^n=B\big(\nabla_x\Phi(\ox,\ow)\R^n+\ker B\big)=B\R^m=\R^s,
$$
which justifies the surjectivity of $\nabla_x\Psi(\ox,\ow)$ and thus completes the proof.$\b$

Recall that the equivalence between ND and SOCQ implies that the first-order qualification condition (\ref{2.9}) automatically holds under ND; see the discussion after the proof of Theorem~\ref{sospwl}. This ensures, by the well-known first-order subdifferential chain rule, that the stationary condition
$\ov\in\partial_x\ph(\ox,\ow)$ via $\ph$ from (\ref{phi}) can be equivalently written as
\begin{equation}\label{ov}
\ov\in\nabla_x\ph_0(\ox,\ow)+\nabla_x\Phi(\ox,\ow)^*\sub\th(\Phi(\ox,\ow)).
\end{equation}
This allows us to consider the KKT system ${\cal P}(w,v)$ given in the form
\begin{eqnarray}\label{kkt}
\left\{\begin{array}{ll}
v=\nabla_xL(x,w,\lm),\quad\lm\in\sub\th(\Phi(x,w))\\
\mbox{with }\;L(x,w,\lm):=\ph_0(x,w)+\langle\lm,\Phi(x,w)\rangle.
\end{array}
\right.
\end{eqnarray}
Similarly, the KKT system for ${\cal P}_r(w,v)$ from {\rm(\ref{fcp3})} is given by
\begin{eqnarray}\label{rkkt}
\left\{\begin{array}{ll}
v=\nabla_x {L}_r(x,w,\mu),\quad\mu\in\sub\vartheta(\Psi(x,w))\\
\mbox{with }\;{L}_r(x,w,\mu):=\ph_0(x,w)+\langle\mu,\Psi(x,w)\rangle.
\end{array}
\right.
\end{eqnarray}

It is easy to see from the reducibility $\th(z)=(\vartheta\circ B)(z)$ around $\oz$ together with the full rank of $B$ that Lagrange multipliers $\lm$ of (\ref{kkt}) and $\mu$ of (\ref{rkkt}) are related by $\lm=B^*\mu$. The next proposition establishes the uniqueness of solutions to (\ref{kkt}) under ND. It is a composite optimization counterpart of \cite[Proposition~4.75]{bs} obtained for optimization problems with constraints $\Phi(x,z)\in\Th$ under the corresponding reducibility and nondegeneracy conditions.

\begin{proposition}\label{uniq}{\bf(uniqueness of Lagrange multipliers for composite problems under ND).} Let $\ox$ be a feasible solution to ${\cal P}(\ow,\ov)$ for the parameter pair $(\ow,\ov)$ with $\ov$ from {\rm (\ref{ov})} and $(\oz,\ov)\in\gph\sub\th$, let {\rm ND} hold, and let $\th\in CPWL$. Then the set of Lagrange multipliers
\begin{equation}\label{lset}
\Big\{\bar\lm\in\sub\th(\Phi(\ox,\ow))\;\,\colon\;\ov=\nabla_x L(\ox,\ow,\bar\lm)\Big\}
\end{equation}
for the KKT system {\rm(\ref{kkt})} is singleton.
\end{proposition}
{\it Proof} Pick $\lm_1,\lm_2$ from (\ref{lset}). It follows from (\ref{lset}) and the subdifferential description (\ref{fos}) for CPWL functions that $\lm_1-\lm_2\in\ker\nabla_x\Phi(\ox,\ow)^*$ and
$$
\lm_s=\sum_{i\in K(\oz)}\eta_{si}a_i+\sum_{i\in I(\oz)}\tau_{si}d_i\;\mbox{ with }\;\sum_{i\in K(\oz)}\eta_{si}=1,\;\;\eta_{si},\tau_{si}\ge 0\;\mbox{ for }\;s=1,2.
$$
Then employing assertions (i) and (ii) of Theorem~\ref{sospwl}, we get
$$
\begin{array}{lll}
\lm_1-\lm_2&=&\disp{\sum_{i\in K(\oz)}\eta_{1i}a_i+\sum_{i\in I(\oz)}\tau_{1i}d_i-\sum_{i\in K(\oz)}\eta_{2i}a_i-\sum_{i\in I(\oz)}\tau_{2i}d_i}\\
&=&\disp{\sum_{j\in K(\oz)}\eta_{2j}\sum_{i\in K(\oz)}\eta_{1i}(a_i-a_j)+\sum_{i\in I(\oz)}\tau_{1i}d_i-\sum_{i\in I(\oz)}\tau_{2i}d_i\in S(\oz)}
\end{array}
$$
thus showing that $\lm_1=\lm_2$ by SOCQ (\ref{2qcf2}), which is equivalent to ND. $\b$

Now we introduce a new second-order condition formulated via the initial data of the composite optimization problem (\ref{fcp2}) and then show that it provides a complete characterization of full stability of local minimizers under ND. This condition is crucial in stability issues for composite optimization playing here the role similar to Robinson's SSOSC \cite{rob} for classical NLPs; therefore, we keep this name in what follows with  adding ``composite."

\begin{definition}{\bf (composite SSOSC).}\label{sssoc} We say that the {\sc composite SSOSC} holds at  $(\ox,\ow,\ov,\bar\lm)\in\R^n\times\R^d\times\R^n\times\R^m$ with $\ov$ and $\bar\lm$ satisfying {\rm(\ref{ov})} and {\rm(\ref{kkt})}, respectively, if
\begin{equation}\label{sssoc1}
\la u,\nabla^2_{xx}L(\ox,\ow,\bar\lm)u\ra>0\;\mbox{ for all }\;0\ne u\in {\cal S},
\end{equation}
where $L$ is the Lagrangian from {\rm(\ref{kkt})} while the subspace ${\cal S}$ is defined by
\begin{equation}\label{sssoc2}
\begin{array}{ll}
{\cal S}:=\Big\{u\in\R^n\,\colon&\la a_i-a_j,\nabla_x\Phi(\ox,\ow)u\ra=0\;\;\mbox{for}\;\;i,j\in\Gamma(J_1),\\
&\la d_t,\nabla_x\Phi(\ox,\ow)u\ra=0\;\;\mbox{for}\;\;t\in\Gamma(J_2)\Big\}
\end{array}
\end{equation}
via the index sets $\Gamma(J_1)$ and $\Gamma(J_2)$ taken from {\rm(\ref{feature})}.
\end{definition}\vspace*{-0.05in}

Observe the following description of the subspace (\ref{sssoc2}) of the positive definiteness of the Lagrangian Hessian in the composite SSOSC:
\begin{eqnarray}\label{compS}
u\in{\cal S}\Longleftrightarrow\nabla_x\Phi(\ox,\ow)u\in\dom\sub^2\th(\ox,\ov),
\end{eqnarray}
which is implied by (\ref{domcod}) and reveals the second-order nature of this subspace. The composite SSOSC reduces to Robinson's SSOSC for NLPs by putting $\Gamma(J_1)=\emp$ and $\Gamma(J_2)=J_2$ in (\ref{sssoc1}), (\ref{sssoc2}). It also reduces to \cite[Definition~6.4]{mrs} and \cite[Definition~7.2]{mrs} in the corresponding settings of MPPCs and ENLPs.

The next lemma is important, together with the second-order subdifferential chain rule, for deriving the aforementioned characterization of full stability of local minimizers in (\ref{fcp2}).

\begin{lemma}\label{zero}{\bf(second-order subdifferential property of CPWL functions).} Take a pair $(\oz,\ov)\in\gph\th$ for a CPWL function $\th$. Then we have $0\in\sub^2\th(\oz,\ov)(u)$ whenever $u\in\dom\sub^2\th(\oz,\ov)$.
\end{lemma}
{\it Proof} Pick $u\in\dom\sub^2\th(\oz,\ov)$ and find $w\in\sub^2\th(\oz,\ov)(u)$, so we deduce from (\ref{2.8}) and (\ref{2nd}) that  $(w,-u)\in N((\oz,\ov),\gph\sub\th)$. Applying (\ref{2nd-val}) gives us a quadruple $(P_1,Q_1,P_2,Q_2)\in{\cal A}$ with
$$
w\in{\cal F}_{\tiny\{P_1,Q_1\},\{P_2,Q_2\}}\;\mbox{ and }\;-u\in{\cal G}_{\tiny\{P_1,Q_1\},\{P_2,Q_2\}}.
$$
Since we always have $0\in {\cal F}_{\tiny\{P_1,Q_1\},\{P_2,Q_2\}}$, it follows that
$$
(0,-u)\in{\cal F}_{\tiny\{P_1,Q_1\},\{P_2,Q_2\}}\times{\cal G}_{\tiny\{P_1,Q_1\},\{P_2,Q_2\}},
$$
which implies by (\ref{2nd-val}) the claimed inclusion $0\in\sub^2\th(\oz,\ov)(u)$. $\b$

We now proceed with establishing the main result of this section, which provides a complete characterization of fully stable local minimizers of  ${\cal P}(\ow,\ov)$  entirely via the initial data.

\begin{theorem}{\bf(second-order characterization of full stability in composite optimization).}\label{ssooc} Let $\ox$ be a feasible solution to ${\cal P}(\ow,\ov)$ from {\rm(\ref{fcp2})} for the parameter pair $(\ow,\ov)$ with $\ov$ from {\rm (\ref{ov})}, let $\th\in CPWL$, and let $(\oz,\ov)\in\gph\sub\th$ with $\oz=\Phi(\ox,\ov)$. Under the validity of ND, let $\bar\lm$ be a unique solution of the KKT system {\rm (\ref{kkt})}. Then $\ox$ is a fully stable local minimizer of ${\cal P}(\ow,\ov)$ if and only if the composite SSOSC from Definition~{\rm\ref{sssoc}} is satisfied.
\end{theorem}
{\it Proof} If $\ox$ is a fully stable local minimizer of ${\cal P}(\ow,\ov)$, then it is also a fully stable local minimizer of the reduced problem ${\cal P}_r(\ow,\ov)$ by Proposition \ref{red1}. It follows from this proposition that the partial Jacobian matrix $\nabla_x\Psi(\ox,\ow)$ of $\Psi=B\circ \Phi$ has full rank. Employing \cite[Theorem~5.1]{mrs} tells us that full stability of $\ox$ for the reduced problem ${\cal P}_r(\ow,\ov)$  is equivalent to
\begin{eqnarray}\label{ss00}
[(p,q)\in{\cal T}_r(\ox,\ow,\ov)(u),\;u\ne0]\Longrightarrow\la p,u\ra>0
\end{eqnarray}
via the set-valued mapping ${\cal T}_r(\ox,\ow,\ov)\colon\R^m\tto\R^{m}\times\R^d$ defined by
\begin{eqnarray*}
\begin{array}{ll}
{\cal T}_r(\ox,\ow,\ov)(u)&\colon=\Big(\nabla^2_{xx}\ph_0(\ox,\ow)u,\nabla^2_{xw}\ph_0(\ox,\ow)u\Big)\\
&+\Big(\nabla^2_{xx}\la\bar\mu,\Psi\ra(\ox,\ow)u,\nabla^2_{xw}\la\bar\mu,\Psi\ra(\ox,\ow)u\Big)\\
&+\Big(\nabla_x\Psi(\ox,\ow),\nabla_w\Psi(\ox,\ow)\Big)^*\partial^2{\vartheta}(\oz,\bar\mu)(\nabla_x\Psi(\ox,\ow)u),\quad u\in\R^m
\end{array}
\end{eqnarray*}
where $\bar\mu$ is a unique solution to the reduced KKT system {\rm (\ref{rkkt})} associated with  $(x,w,v):=(\ox,\ow,\ov)$. The full rank of $\nabla_x\Psi(\ox,\ow)$ allows us to use the second-order chain rule from \cite[Theorem~3.1]{mr} and get
\begin{eqnarray*}
{\cal T}_r(\ox,\ow,\ov)(u)=\Big(\nabla^2_{xx}\ph_0(\ox,\ow)u,\nabla^2_{xw}\ph_0(\ox,\ow)u\Big)+D^*{\partial}_x(\vartheta\circ\Psi)(\ox,\ow,\ov)(u).
\end{eqnarray*}
By the representation $(\vartheta\circ\Psi)(x,w)=(\th\circ\Phi)(x,w)$ around $(\ox,\ow)$ we have
\begin{equation}\label{th3.1}
{\cal T}_r(\ox,\ow,\ov)(u)=\Big(\nabla^2_{xx}\ph_0(\ox,\ow)u,\nabla^2_{xw}\ph_0(\ox,\ow)u\Big)+D^*{\partial}_x(\th\circ\Phi)(\ox,\ow,\ov)(u).
\end{equation}
Applying the second-order chain rule from Corollary~\ref{soc} to $\th\circ\Phi$ in (\ref{th3.1}) together  with (\ref{ss00}) tells us that  $\ox$ being a fully stable local minimizer of the reduced problem ${\cal P}_r(\ow,\ov)$ is equivalent to the validity of the inequality
\begin{equation}\label{ss01}
\la u,\nabla^2_{xx}L(\ox,\ow,\bar\lm)u\ra+\la q,\nabla_x\Phi(\ox,\ow)u\ra>0\;\mbox{if}\;q\in\partial^2\th(\oz,\bar\lm)(\nabla_x\Phi(\ox,\ow)u),u\ne 0.
\end{equation}
Pick $0\ne u\in {\cal S}$ and get by (\ref{compS}) that $\nabla_x\Phi(\ox,\ow)u\in \dom\sub^2\th(\oz,\ov)$. Thus it follows from Lemma~\ref{zero} that $0\in\partial^2\th(\oz,\bar\lm)(\nabla_x\Phi(\ox,\ow)u)$ implying by (\ref{ss01}) that
$$
\la u,\nabla^2_{xx}L(\ox,\ow,\bar\lm)u\ra=\la u,\nabla^2_{xx}L(\ox,\ow,\bar\lm)u\ra+\la 0,\nabla_x\Phi(\ox,\ow)u\ra>0
$$
verifying therefore the ``only if" statement.

To establish next the ``if" part of the theorem, assume that  $u\ne 0$ and that $q\in\partial^2\th(\oz,\bar\lm)(\nabla_x\Phi(\ox,\ow)u)$, which yields $u\in{\cal S}$. Then \cite[Theorem~2.1]{pr98} together with the convexity of $\th$ ensures that $\la q,\nabla_x\Phi(\ox,\ow)u\ra\ge 0$, and hence we have
$$
\la u,\nabla^2_{xx}L(\ox,\ow,\bar\lm)u\ra+\la q,\nabla_x\Phi(\ox,\ow)u\ra\ge\la u,\nabla^2_{xx}L(\ox,\ow,\bar\lm)u\ra>0
$$
by the assumed composite SSOSC. This implies by (\ref{ss01}) that $\ox$ is a fully stable local minimizer of the reduced problem ${\cal P}_r(\ow,\ov)$. Appealing finally to Proposition~\ref{red1} shows that $\ox$ is a fully stable local minimizer of problem ${\cal P}(\ow,\ov)$ and thus completes the proof of the theorem. $\b$

The obtained characterization extends the results of \cite[Theorem~6.6]{mrs} for MPPCs, of \cite[Theorem~7.3]{mrs} for ENLPs, and of \cite[Theorem~6.3]{ms14} for unconstrained minimax problems. An important advantage of Theorem~\ref{ssooc} is that it allows us to characterize full stability of local minimizers in (nonsmooth) minimax problems with polyhedral constraints, which is done in the next section while cannot be obtained by using the developments in \cite{mrs,ms14}.\vspace*{-0.2in}

\section{Full Stability in Constrained Minimax Problems}\label{sec:minmax}

This section deals with applications of Theorem~\ref{ssooc} and second-order subdifferential calculations from \cite{ms15} to characterizing fully stable local minimizers for the following class of {\em minimax} problems with {\em polyhedral constraints}:
\begin{equation}\label{minmax}
\mbox{minimize }\;\max\{\ph_1(x),\ldots,\ph_l(x)\}\;\mbox{ over }\;\Upsilon(x):=(\zeta_1(x),\ldots,\zeta_r(x))\in Z
\end{equation}
with $r+l=m$, where the functions $\ph_i\colon\R^n\to\R$ for $i=1,\ldots,l$ and $\zeta_s\colon\R^n\to\R$ for $s=1,\ldots,r$ are ${\cal C}^2$-smooth around the reference points,
and where the convex polyhedron $Z\subset\R^r$ is given by
\begin{equation}\label{zpol}
Z:=\Big\{y\in\R^r\,\colon\;\la c_t,y\ra\le\tau_t\;\;\mbox{for all}\;\;t=1,\ldots,p\;\Big\}
\end{equation}
with $(c_t,\tau_t)\in\R^r\times\R$ for $t=1,\ldots,p$. The minimax counterpart of ${\cal P}(w,v)$ from above is written as: minimize
\begin{equation}\label{minmax2}
\max\{\ph_1(x,w),\ldots,\ph_l(x,w)\}+\dd(\Upsilon(x,w);Z)-\la v,x\ra\;\mbox{s.t.}\;x\in\R^n
\end{equation}
with $(w,v)\in\R^d\times\R^n$. We say that $x\in\R^n$ is a feasible point to it (\ref{minmax2}) if $\Upsilon(x,w)\in Z$. Note that problem (\ref{minmax2}) differs from  ${\cal P}(w,v)$ in (\ref{fcp2}) due to {\em nonsmoothness} of all the summands in  (\ref{minmax2}) but $\la v,x\ra$. Let us show that nevertheless (\ref{minmax2}) can be reduced to the composite form (\ref{fcp2}) as follows. Consider the mapping $\Phi\colon \R^n\times\R^d\to\R^{l+r}=\R^m$ given by
\begin{equation}\label{ss04}
\Phi(x,w):=(\Xi(x,w),\Upsilon(x,w))\;\mbox{ for all }\;(x,w)\in\R^n\times\R^d
\end{equation}
with the mapping $\Upsilon$ taken from (\ref{minmax}) and $\Xi(x,w):=(\ph_1(x,w),\ldots,\ph_l(x,w))$. Remembering that $r+l=m$, define $\th\colon\R^{l+r}\to\Bar\R$ by
\begin{equation}\label{theta2}
\left\{\begin{array}{ll}
\th(x):=\max\Big\{\la a_1,x\ra,\ldots,\la a_l,x\ra\Big\}+\dd(x;{\cal Z})\;\mbox{ for }\;x\in\R^{l+r}=\R^m\\
\mbox{with}\quad{\cal Z}:=\Big\{x\in\R^{l+r}\,\colon\;\la d_t,x\ra\le\tau_t\;\mbox{ for }\;t=1,\ldots,p\;\Big\},
\end{array}
\right.
\end{equation}
where the generating vectors $a_i$ and $d_t$  are constructed from the unit vectors $e_i\in\R^l$ and the vectors $c_t\in\R^r$ from (\ref{zpol}) by, respectively,  \begin{equation}\label{ss05}
a_i:=(e_i,0)\;\mbox{ for }\;i=1,\ldots,l\;\mbox{ and }\;d_t:=(0,c_t)\;\mbox{ for }\;t=1,\ldots,p.
\end{equation}
Observe the $\th$ from (\ref{theta2}) is a CPWL function in the summation form (\ref{theta}). Thus we can represent the constrained minimax problem (\ref{minmax2}) in the composite optimization form (\ref{fcp2}) written as
\begin{equation}\label{minmax3}
\mbox{minimize }\;(\th\circ\Phi)(x,w)-\la v,x\ra\;\mbox{ subject to }\;x\in\R^n
\end{equation}
with $\th$ taken from (\ref{theta2}) with parameters (\ref{ss05}) and the ${\cal C}^2$-smooth mapping (\ref{ss04}).

Now we can apply Theorem~\ref{ssooc} to (\ref{minmax3}) and derive in this way a second-order characterization of full stability of local solutions to the minimax problem (\ref{minmax2}) via its initial data. Prior to that, let us specify the nondegeneracy condition ND for problem (\ref{minmax3}) and presents it in terms of the original minimax problem (\ref{minmax2}) without appealing to the matrix $B$ from the proof of Lemma~\ref{fred}.

Denote $\oz_1:=\Xi(\ox,\ow)$ and $\oz_2:=\Upsilon(\ox,\ow)\in Z$ and construct the index sets
\begin{equation}\label{active3}
\begin{array}{lll}
{\cal K}(\oz_1):&=&\Big\{i\in\{1,\ldots,l\}\,\colon\;\max\{\ph_1(\ox,\ow),\ldots,\ph_l(\ox,\ow)\}=\ph_i(\ox,\ow)\Big\},\\
{\cal I}(\oz_2):&=&\Big\{t\in\{1,\ldots,p\}\,\colon\;\la c_t,\oz_2\ra=\tau_t\Big\}
\end{array}
\end{equation}
via the data of (\ref{minmax}) and (\ref{zpol}). It is easy to observe that ${\cal K}(\oz_1)=K(\oz)$ and ${\cal I}(\oz_2)=I(\oz)$ for the index sets defined in (\ref{active2}) for the function $\th$ from (\ref{theta2}) with $\oz:=(\oz_1,\oz_2)\in\dom\th$.

\begin{proposition}{\bf (equivalent form of qualification condition ND for constrained minimax problems).}\label{ndmin0} Let $\ox$ be a feasible solution to {\rm(\ref{minmax2})} corresponding to $(\ow,\ov)$, and let $\oz=(\oz_1,\oz_2)$ with $\oz_1=\Xi(\ox,\ow)$, $\oz_2=\Upsilon(\ox,\ow)$, and $(\oz,\ov)\in\gph\sub\th$, where the mappings $\Xi$ and $\Upsilon$ and the CPWL function $\th$ are defined by {\rm(\ref{ss04})} and {\rm(\ref{theta2})}, respectively. Then the nondegeneracy condition ND in the framework of the minimax problem {\rm (\ref{minmax2})} can be equivalently written as
\begin{equation}\label{ndmin}
{\cal D}\cap\ker\Big(\nabla_x\Xi(\ox,\ow)^*,\nabla_x\Upsilon(\ox,\ow)^*\Big)=\{0\}
\end{equation}
with the set ${\cal D}$ given by
$$
\begin{array}{lll}
{\cal D}:=\Big\{(y_1,y_2)\in\R^l\times\R^r\;\colon &y_1\in\span\{e_i-e_j\;\colon\;i,j\in{\cal K}(\oz_1)\},\;\mbox{and}\\
&y_2\in\span\{c_t\;\colon\;t\in{\cal I}(\oz_2)\}\Big\},
\end{array}
$$
 where $e_i\in\R^l$ are the unit vectors, and where the index sets ${\cal K}(\oz_1)$ and ${\cal I}(\oz_2)$ are defined in {\rm(\ref{active3})}.
\end{proposition}
{\it Proof} Applying the nondegeneracy condition ND to the composite optimization form (\ref{minmax3}) of the minimax problem  (\ref{minmax3}) and using Theorem~\ref{fred2} on the equivalence of ND to SOCQ give us
$$
\sub^2\th(\oz,\ov)(0)\cap\ker\Phi(\ox,\ow)^*=\{0\}
$$
with $\Phi$ and $\th$ taken from {\rm(\ref{ss04})} and {\rm(\ref{theta2})}, respectively. Then the second-order calculations of Theorem~\ref{sospwl} together with the equalities ${\cal K}(\oz_1)=K(\oz)$ and ${\cal I}(\oz_2)=I(\oz)$ reveal that
\begin{equation}\label{ss08}
\begin{array}{lll}
\sub^2\th(\oz,\ov)(0)&=\span\Big\{a_i-a_j\,\colon\;i,j\in K(\oz)\Big\}+\span\Big\{d_t\,\colon\;t\in I(\oz)\Big\}\\
&=\Big\{(y_1,y_2)\in\R^l\times\R^n\,\colon \;y_1\in\span\{e_i-e_j\;\colon\;i,j\in{\cal K}(\oz_1)\},\\
&\hspace{3.7cm}y_2\in\span\{c_t\;\colon\;t\in{\cal I}(\oz_2)\}\Big\},
\end{array}
\end{equation}
where the vectors $a_i$ and $d_t$ are taken from (\ref{ss05}). Observing that
$$
\nabla_x\Phi(\ox,\ow)=\left(\begin{array}{c}
\nabla_x\Xi(\ox,\ow)\\
\nabla_x\Upsilon(\ox,\ow)
\end{array}
\right)
$$
and combining this with representation (\ref{ss08}) justify the equivalent form (\ref{ndmin}) of the ND condition in the minimax problem under consideration. $\b$

After these adjustments, we now derive a characterization of fully stable local minimizers of (\ref{minmax2}). The KKT system for (\ref{minmax2}) can be expressed as
\begin{equation}\label{kktmin}
\left\{\begin{array}{ll}
\disp{\ov=\sum_{i=1}^{l}\bar\lm_i\nabla_x\ph_i(\ox,\ow)+\sum_{s=1}^{r}\bar\mu_s\nabla_x\zeta_s(\ox,\ow)}\\
\mbox{with}\quad\disp{\bar\lm_i\ge 0,\;\sum_{i=1}^{l}\bar\lm_i=1,\;(\bar\mu_1,\ldots,\bar\mu_r)\in N(\oz_2;Z)},
\end{array}\right.
\end{equation}
where $\oz_2=\Upsilon(\ox,\ow)$, and where $Z$ is taken from (\ref{zpol}). The following definition is an adaptation of the composite SSOSC for the minimax problem (\ref{minmax2}).

\begin{definition}{\bf (minimax SSOSC).}\label{sssocm} Given $\varpi:=(\bar\lm,\bar\mu)\in\R^{l}\times\R^r$ and $\ov$ from {\rm(\ref{kktmin})}, we say that the
{\sc minimax SSOSC} holds at $(\ox,\ow,\ov,\varpi)$ if
\begin{equation}\label{sssocm1}
\sum_{i=1}^{l}\bar\lm_i\la u,\nabla^2_{xx}\ph_i(\ox,\ow)u\ra+\sum_{s=1}^{r}\bar\mu_s\la u,\nabla^2_{xx}\zeta_s(\ox,\ow)u\ra>0\;\mbox{ for all }\;0\ne u\in{\cal S},
\end{equation}
where the subspace ${\cal S}$ is defined by
\begin{eqnarray*}
\begin{array}{ll}
{\cal S}:=\Big \{u\in\R^n\,\colon&\la\nabla_x\ph_i(\ox,\ow),u\ra=\gg\;\;\mbox{for}\;\;i\in\Gamma(J_1)\;\;\mbox{and}\\
&\la d_t,\nabla_x\Upsilon(\ox,\ow)u\ra=0\;\;\mbox{for}\;\;t\in\Gamma(J_2)\Big\}
\end{array}
\end{eqnarray*}
via the index sets $\Gamma(J_1)$ and $\Gamma(J_2)$ taken from {\rm(\ref{feature})} and some constant $\gg\in\R$.
\end{definition}

Next we extend \cite[Theorem~6.3]{mr} to {\em constrained} minimax problems.

\begin{theorem}{\bf(characterization of fully stable solutions to constrained minimax problems).}\label{fulmin} Let $\ox$ be a feasible solution to the minimax problem {\rm(\ref{minmax2})} corresponding to $(\ow,\ov)$ with $\ov\in\sub_x(\th\circ\Phi)(\ox,\ow)$, where $\th$ and $\Phi$ are taken from {\rm(\ref{theta2})} and {\rm(\ref{ss04})}, respectively. Assume that the ND condition {\rm(\ref{ndmin})} holds, and let $\varpi=(\bar\lm,\bar\mu)\in\R^{l}\times\R^r$ be a unique solution to {\rm (\ref{kktmin})}. Then $\ox$ is a fully stable local minimizer of {\rm(\ref{minmax2})} if and only if the minimax SSOSC from {\rm(\ref{sssocm1})} holds.
\end{theorem}
{\it Proof} Following the lines above, we implement Theorem~\ref{ssooc} in the constrained minimax setting by observing that the Lagrangian for (\ref{minmax2}) can be represented as $L(\ox,\ow,\varpi)=\la\varpi,\Phi(\ox,\ow)\ra$. This gives us
$$
\nabla^2_{xx}L(\ox,\ow,\varpi)=\sum_{i=1}^{l}\bar\lm_i\nabla^2_{xx}\ph_i(\ox,\ow)+\sum_{s=1}^{r}\bar\mu_s\nabla^2_{xx}\zeta_s(\ox,\ow).
$$
Furthermore, it is easy to see that the set $\cal S$ from Definition~\ref{sssocm} is an adaptation of $\cal S$ from (\ref{sssoc2}) to (\ref{minmax2}). Thus the claimed second-order characterization of full stability  in (\ref{minmax2}) readily follows from the equivalence in Theorem~\ref{ssooc}. $\b$\vspace*{-0.2in}

\section{Conclusions}
This paper provides various applications of the second-order subdifferential theory for CPWL functions recently developed in \cite{ms15}. The obtained results prove the importance and power of such constructions in variational analysis. Following this way, we plan to proceed with further applications. In particular, our intention is establish the equivalence between the Lipschitz-like/Aubin property  and Robinson's strong regularity in the CPWL framework.\vspace*{-0.1in}

\begin{acknowledgements}
The authors are grateful to the referees for their helpful remarks.
This research was partly supported by the National Science Foundation under grants DMS-1007132 and DMS-1512846 and by the Air Force Office of Scientific Research grant \#15RT0462.\vspace*{-0.1in}
\end{acknowledgements}

\begin{thebibliography}{}

\bibitem{dem} Demyanov, V.F., Rubinov, A.M.: Constructive Nonsmooth Analysis. Peter Lang, Frankfurt
(1995)

\bibitem{rw} Rockafellar, R.T., Wets, R.J-B.: Variational Analysis. Springer, Berlin (1998)

\bibitem{ms15} Mordukhovich, B.S., Sarabi, M.E.: Generalized differentiation of piecewise linear functions
in second-order variational analysis. Nonlinear Anal. 132, 240--273 (2016)

\bibitem{m92} Mordukhovich, B.S.: Sensitivity analysis in nonsmooth optimization. SIAM Proc. Appl.
Math. 58, 32--46 (1992)

\bibitem{mr} Mordukhovich, B.S., Rockafellar, R.T.: Second-order subdifferential calculus with application
to tilt stability in optimization. SIAM J. Optim. 22, 953--986 (2012)

\bibitem{mnn} Mordukhovich, B.S., Nam, N.M., Nhi, N.T.Y.: Partial second-order subdifferentials in
variational analysis and optimization. Numer. Func. Anal. Optim. 35, 1113--1151 (2014)

\bibitem{lpr} Levy, A.B., Poliquin, R.A., Rockafellar, R.T.: Stability of locally optimal solutions.
SIAM J. Optim. 10, 580--604 (2000)

\bibitem{m06} Mordukhovich, B.S.: Variational Analysis and Generalized Differentiation, I: Basic Theory;
II: Applications. Springer, Berlin (2006)

\bibitem{rob84} Robinson, S.M.: Local structure of feasible sets in nonlinear programming, ii: Nondegeneracy.
Math. Program. Stud. 22, 217--230 (1984)

\bibitem{mrs} Mordukhovich, B.S., Rockafellar, R.T., Sarabi, M.E.: Characterizations of full stability
in constrained optimization. SIAM J. Optim. 23, 1810--1849 (2013)

\bibitem{bs} Bonnans, J.F., Shapiro, A.: Perturbation Analysis of Optimization Problems. Springer,
New York (2000)

\bibitem{mos} Mordukhovich, B.S., Outrata, J.V., Sarabi, M.E.: Full stability of locally optimal solution
in second-order cone programming. SIAM J. Optim. 24, 1581--1613 (2014)

\bibitem{r} Rockafellar, R.T.: Extended nonlinear programming. Nonlinear Optimization and Related
Topics 36, 381--399 (2000)

\bibitem{mn3} Mordukhovich, B.S., Nghia, T.T.A.: Full Lipschitzian and Holderian stability in optimization
with applications to mathematical programming and optimal control. SIAM
J. Optim. 24, 1344--1381 (2014)

\bibitem{mnr} Mordukhovich, B.S., Nghia, T.T.A., Rockafellar, R.T.: Full stability in finitedimensional
optimization. Math. Oper. Res. 40, 226--252 (2015)

\bibitem{ms14} Mordukhovich, B.S., Sarabi, M.E.: Variational analysis and full stability of optimal
solutions to constrained and minimax problems. Nonlinear Anal. 121, 36--53 (2015)

\bibitem{mn14} Mordukhovich, B.S., Nghia, T.T.A.: Local strong maximal monotonicity and full
stability for parametric variational systems. to appear in SIAM J. Optim.,
http://www.optimization-online.org/DB-HTML/2015/09/5089.html

\bibitem{rob} Robinson, S.M.: Strongly regular generalized equations. Math. Oper. Res. 5, 43--62
(1980)

\bibitem{pr98} Poliquin, R.A., Rockafellar, R.T.: Tilt stability of a local minimum. SIAM J. Optim. 8,
287--299 (1998)






\end{thebibliography}

\end{document}